\documentclass{amsart}
\usepackage{amssymb}
\usepackage{amsmath}
\usepackage{amscd}
\usepackage{graphicx}
\usepackage{latexsym}
\usepackage{amsrefs}
\usepackage{hyperref}
\usepackage{setspace}

\def\Z{\mathbb{Z}}

\def\M{{\mathcal M}}

\newtheorem{theorem}{Theorem}[section]
\newtheorem{lemma}[theorem]{Lemma}

\newtheorem{corollary}[theorem]{Corollary}

\newtheorem{remark}[theorem]{Remark}
\newtheorem{example}[theorem]{Example}

\def\CP{\hbox{${\mathbb C}P^2$}}
\def\CPB{\hbox{$\overline{{\mathbb C}P^2}$}}
\def\CPthree{\hbox{${\mathbb C}P^3$}}

\def\w{\omega}

\title{Count of Genus Zero J-Holomorphic Curves in Dimensions Four and Six}
\author{Ahmet Beyaz}
\address{Department of Mathematics, Middle East Technical University, Ankara 06800 Turkey}
\email{beyaz@metu.edu.tr}
\subjclass[2010]{53D05, 57R55, 53D45}
\keywords{symplectic $6$-manifolds, $J$-holomorphic curves}

\begin{document}
\begin{abstract}  In this note, genus zero Gromov-Witten invariants are reviewed and then applied in some examples of dimension four and six. It is also proved that the use of genus zero Gromov-Witten invariants in the class of embedded $J$-holomorphic curves to distinguish the deformation types of symplectic structures on a smooth $6$-manifold is restricted in the sense that they can not distinguish the symplectic structures on $X_1\times S^2$ and $X_2\times S^2$ for two minimal, simply connected, symplectic $4$-manifolds $X_1$ and $X_2$ with $b_2^+(X_1)>1$ and $b_2^+(X_2)>1$.
\end{abstract}
\maketitle
\setcounter{section}{-1}

\section{Introduction}

The count of $J$-holomorphic curves in a symplectic manifold carries information about the properties of the symplectic structure on the manifold. This kind of study was first established by M. Gromov (\cite{Gromov1985}) in 1985. Later it was improved and applied in many ways (\cite{Kontsevich1994, Ruan1995, Ruan1999, Fukaya1999, Zinger2008}), in particular for symplectic $6$-manifolds. 

The study of symplectic structures in dimension six is so incomplete that even for manifolds with simple topology like the homotopy projective spaces, it is not clear if they admit any symplectic structures except $\CPthree$ itself. Nevertheless it is conjectured that given a topological $4$-manifold $X$ which admits symplectic structures, the classification of smooth structures on $X$ is equivalent to the classification of deformation types of symplectic structures on $X\times S^2$ (\cite{McDuff1998} page 437, \cite{Ruan2002}). It is shown that this conjecture holds for elliptic surfaces by Ruan and Tian (\cite{Ruan1997}). Moreover it was proven that if $X_1\times S^2$ and $X_2\times S^2$ are deformation equivalent, then some branched covers of $X_1$ and $X_2$ are diffeomorphic (\cite{Ruan2002}).

In Section~\ref{GW} genus zero Gromov-Witten invariants of a symplectic manifold are defined. In the subsequent sections, we explain how these invariants are applied in some examples. The invariants are defined using simple (i.e. non multiply covered) curves. Taking multiply covered curves into account does not change the results of this paper regarding minimal $4$-manifolds.

In Remark~\ref{SW} the relation of Gromov-Witten invariants of a $4$-manifold $X$ to its Seiberg-Witten invariants is discussed in a nutshell. This suggests that some of Gromov-Witten invariants of $X\times S^2$ can be given in terms of the Seiberg-Witten invariants of $X$. However our results imply that in genus zero case one can not get information except for the class of an exceptional sphere. In particular, for a minimal symplectic $4$-manifold, Seiberg-Witten invariants do not contribute to genus zero Gromov-Witten theory. In the last section we prove Theorem~\ref{main}, which shows that the efficiency of genus zero Gromov-Witten invariants to distinguish the symplectic deformation types on a smooth $6$-manifold is restricted in the sense that they can not distinguish the symplectic structures on $X_1\times S^2$ and $X_2\times S^2$ for two minimal, simply connected, symplectic $4$-manifolds $X_1$ and $X_2$, with $b_2^+(X_1)>1$ and $b_2^+(X_2)>1$ in the pushforwards of second homology classes with minimal genus zero. This is a consequence of Theorem~\ref{main}. Another consequence is Corollary~\ref{deformationEq} which states that if $X_1$ and $X_2$ are homeomorphic  symplectic $4$-manifolds with $b_2^+(X_1)>1$ and $b_2^+(X_2)>1$ and if $X_1$ is not minimal and $X_2$ is minimal, then $X_1\times S^2$ and $X_2\times S^2$ are not symplectic deformation equivalent.

In the following discussion, as a convention, $g(A)$ will be the genus of a surface and $g([A])$ will be the minimal genus of embedded representatives of $[A]\in H^2(X;\Z)$ in a $4$-manifold $X$. The self intersection of $A$ and $[A]$ in $X$ is denoted by $A^2$ and $[A]^2$, respectively.

\section{Genus Zero Gromov-Witten Invariants} \label{GW}

In this section the definition of Gromov-Witten invariants in genus zero is reviewed. A Gromov-Witten type invariant is an invariant of the deformation type of symplectic structures on a manifold. Gromov-Witten invariants of symplectic manifolds count the number of connected $J$-holomorphic curves in a particular homology class which pass through a number of points. In practice, for genus zero invariants, for a given homology class this count is done by tracing the oriented intersection points of a moduli space and a number of cohomology elements. A compact symplectic manifold $(M,\w)$ is called semipositive if there are no spherical homology classes $[A]\in H_2(M;\Z)$ such that $\w([A])>0$ and $2-n<c_1(M)[A]<0$. In particular if the dimension of the manifold in these definitions is less than or equal to six, then the manifold is semipositive. 

\subsection*{Definitions} 

For a symplectic manifold $(M,\omega)$ (or just $M$ when there is no ambiguity), let $J$ be a generic compatible almost complex structure. A $J$-holomorphic curve in $M$ is a smooth map $u$ from a genus $g$ complex curve into $M$ such that $J \circ du= du \circ i$ where $i$ is the complex structure on the curve. 

Given a nonzero homology class $[A] \in H_2(M;\Z)$ and a positive integer $k$, consider the moduli space $\M^M_{[A],g,k}$ of all simple genus $g$ maps into $M$ with $k$ distinct marked points where the homology class of the image is $[A]$, up to reparametrization. $\M^M_{[A],g,k}$ consists of the equivalence classes $[u, x_1, \cdots, x_k]$. Since we deal with genus zero Gromov-Witten invariants, i.e. $g=0$, we drop the subscript $g$ from the notation. 

Let $M^k$ denote the Cartesian product of $k$ copies of $M$ for $k>0$. For $1<j<k$, let $\pi_j$ denote the projection map $M^k\rightarrow M$ onto the $j$th factor. The evaluation map $ev:{\M}^M_{[A],k} \hookrightarrow M^k$ is defined by $ev([u,x_1,\cdots, x_k])=(u(x_1),\cdots,u(x_k))$. 

Theorem~6.6.1 of \cite{McDuff2004} states that if $M$ is semipositive and if $[A]\in H_2(M;\Z)$ satisfying Condition \ref{multipleTori} below, then the evaluation map $ev:\M^M_{[A],k} \hookrightarrow M^k$ is a pseudocycle of real dimension $2n-6+2c_1(M)[A]+2k$ where $2n$ is the dimension of the manifold $M$. For the definition of a pseudocycle see Definition~6.5.1 of \cite{McDuff2004}.

\begin{equation}\label{multipleTori} [A]=m[B] \textrm{ and } c_1(M)[B]=0 \Rightarrow m=1
\end{equation}
\noindent for all $m>0$ and for all spherical homology classes $[B]\in H_2(M;\Z)$.

Condition \ref{multipleTori} on $[A]$ is not restrictive in the context of this paper because such classes can not be represented by embedded $J$-holomorphic spheres in a symplectic $4$-manifold with $b_2^+>1$.

According to Theorem~7.1.1 of \cite{McDuff2004}, the $k$-pointed genus zero Gromov-Witten invariant of $(M,\omega)$ in the class $[A]$ is defined as

$$GW^M_{[A],k}(\alpha_1,\cdots,\alpha_k)= ev\cdot f$$

\vspace{2mm}

\noindent where $\alpha_1,\cdots, \alpha_k$ are cohomology classes of $M$ and $f$ is a pseudocycle which is Poincar\'e dual to $\pi_1^*(\alpha_1)\cup\cdots\cup \pi_k^*(\alpha_k)$. See Lemma~6.5.5 of  (\cite{McDuff2004}).

To get nonzero invariants, the sum of the degrees of the cohomology elements must be equal to the dimension of ${\M}^M_{[A],k}$, which is $2n-6+2c_1(M)[A]+2k$. This is called the dimension condition. 

Gromov-Witten invariants can be consistently extended to the case where $k$ is zero.  If $[A]$ is zero, $GW^M_{0,0}$ is set as zero. When $k$ is zero, $M^k$ is a point and any pseudocycle is trivial. If $[A]$ is nonzero, for the dimension condition to be satisfied $2n-6+2c_1(M)[A]$ should be zero. When the dimension $2n$ of $M$ is four, under the assumption that $[A]$ has an embedded sphere representative, this implies that $[A]^2$ should be $-1$ and $GW^M_{[A],0}$ counts the exceptional spheres in the class $[A]$. This is either zero or one (not $-1$ as a convention) as in Example~\ref{example4}. If $2n$ is six, then for $GW^M_{[A],0}$ to be nonzero, $c_1(M)[A]$ should be zero.

Two facts about Gromov-Witten invariants which are used in the subsequent sections are the following two lemmas which are known as the fundamental class axiom and the divisor axiom for Gromov-Witten invariants.

{\lemma \label{fundamentalClassAxiom} Let $M$ be a semipositive symplectic manifold, $[A]$ be a nonzero second homology class and $k\geq1$. Then $GW_{[A],k}(\alpha_1,\cdots,\alpha_{k-1},PD([M]))$ is zero. In other words, there can not be a degree zero cohomology class among $\alpha_i$'s.}

{\lemma \label{divisorAxiom} Let $M$ be a semipositive symplectic manifold, $[A]$ be a nonzero second homology class and $k\geq1$. If the degree of $\alpha_k$ is two, then 
$$GW_{[A],k}(\alpha_1,\cdots,\alpha_{k-1},\alpha_k))=(\alpha_k\cdot[A])\,\,GW_{[A],{k-1}}(\alpha_1,\cdots,\alpha_{k-1})$$}

\section{Dimension Four} In the proofs, we are going to apply different results which appeal to a generic set of compatible almost complex structures. These are Baire sets as well as their intersections, thus one can find a compatible almost complex structure $J$ which is in all of these sets (\cite{Wendl} page 109). 

\subsection*{Exceptional Spheres}
We start this subsection with an example on calculations of Gromov-Witten invariants.

{\example \label{example4} Let $Y$ be a simply connected symplectic $4$-manifold and $X$ be its blowup. Topologically $X$ is diffeomorphic to $Y\#\CPB$. In the blowup of a $4$-manifold, there is an exceptional sphere which is a smooth sphere with self intersection $-1$. Let's find the Gromov-Witten invariant of $X$ for the homology class $[E]$ of the exceptional sphere in $H_2(X;\Z)$ with no other constraints. This means the number of points through which it passes is zero. By the adjunction formula for symplectic $4$-manifolds, $c_1(X)[E]$ is calculated as one. The expected dimension of the moduli space is $2n-6+2c_1(X)[E]$ which is equal to zero. The moduli space ${\M}^X_{[E],0}$ is a finite set of points with orientation. By the positivity of intersections of $J$-holomorphic curves in an almost complex manifold, there is only one $J$-holomorphic curve $E$ which represents $[E]$ in $X$. So $GW_{[E],0}^X$ is $1$. The next theorem is an extension of this example to positive values of $k$. Compare \cite{Ruan2002}.}

{\theorem \label{pm1} Let $Y$ be a symplectic $4$-manifold and $X$ be $Y\#\CPB$. If $[E]$ is the homology class of the exceptional sphere, then the Gromov-Witten invariant $GW^{X}_{[E],k}(PD[E],\cdots,PD[E])$ is equal to $(-1)^k$.}
\begin{proof} Let $J$ be a compatible almost complex structure on $X$. By positivity of intersections of $J$-holomorphic curves in an almost complex manifold, there is only one $J$-holomorphic curve which represents $[E]$ in $X$, which will be denoted by $E$. $GW_{[E],0}^X$ is one. $PD([E])\cdot[E]$ is $-1$. Applying the divisor axiom (\ref{divisorAxiom}) inductively, we find $GW^{X}_{[E],k}(PD[E],\cdots,PD[E])=(-1)^k$.
\end{proof}

{\theorem Let $k>0$. If $Y$ is a simply connected, symplectic $4$-manifold and $X$ is $Y\#\CPB$, then the Gromov-Witten invariant $GW^{X}_{[E],k}(\alpha_1,\cdots,\alpha_k)$ is zero unless $\alpha_i\in H^2(X;\Z)$. If $\alpha_i\in H^2(X;\Z)$ for all $i\in\Z$ such that $1\leq i\leq k$ ($k>0$), then 
$$GW^{X}_{[E],k}(\alpha_1,\cdots,\alpha_k) = (\alpha_1\cdot [E])\cdots(\alpha_k\cdot [E])$$}
\begin{proof} By the dimension condition, the sum of degrees of $\alpha_i$ should be equal to the dimension of the moduli space ${\M}^{X}_{[E],k}$, which is $2n-6+2c_1(X)[E]+2k=2k$. There is no odd degree cohomology because $X$ is simply connected. By the fundamental class axiom (Lemma~\ref{fundamentalClassAxiom}), in order to get a nonzero invariant, all classes must be of degree two. By the divisor axiom the result follows.
\end{proof}

\subsection*{Nonzero Invariants}

This subsection is on the conditions for the invariants to be nonzero. The following lemma is critical in the proofs of the main theorems.

{\lemma \label{A2<0} Let $X$ be a symplectic $4$-manifold with $b_2^+(X)>1$ and $J$ be a generic almost complex structure on $X$ which is compatible with the symplectic structure. Let $A\subset X$ be a connected $J$-holomorphic sphere in the class of $[A]\in H_2(X;\Z)$ such that $[A]$ is nonzero and $[A]$ can be represented by an embedded, connected $J$-holomorphic sphere. Then $A^2$ is less than or equal to $-1$. Moreover $A$ is embedded and it is multiple cover of an exceptional sphere.}
\begin{proof} If $A^2\geq0$, then by the assumption of the lemma there is an embedded $J$-holomorphic sphere in the class $[A]$ with self intersection greater than or eqaul to zero. Since a homologically essential embedded $2$-sphere in a symplectic $4$-manifold with $b^+>1$ always has negative self intersection (\cite{McDuff1990a}), we necessarily have $A^2<0$. If $A^2$ is less than $-1$, $A$ is a multiple cover of an exceptional sphere by Theorem~$1.2$ of \cite{McDuff1996} and it is embedded.
\end{proof}

The next theorem is one of the main results.

{\theorem \label{4A=E} Let $X$ be a symplectic $4$-manifold with $b_2^+(X)>1$ and $J$ be a generic almost complex structure on $X$ which is compatible with the symplectic structure. Let $[A]$ be a  nonzero second homology class of $X$ which can be represented by an embedded, connected $J$-holomorphic sphere. Let $\alpha_1,\cdots,\alpha_k$ be cohomology classes of $X$. If $GW^{X}_{[A],k}(\alpha_1,\cdots,\alpha_k)$ is nonzero, then $[A]$ is the class of an exceptional sphere in $X$.}
\begin{proof} If the Gromov-Witten invariant is nonzero then $[A]$ must have a connected $J$-holomorphic representative, say $A$. 
Lemma~\ref{A2<0} implies that $A=mE$ ($m>1$) for some exceptional sphere $E$. Since only simple curves are considered $m$ should be one and therefore $[A]$ is the class of an exceptional sphere.
\end{proof}

The case where $[A]$ is the zero class is excluded. We refer the reader to a more general source (\cite{McDuff2004}) for a discussion on the zero class.

{\remark \label{SW} These results are compatible with the results of Taubes (\cite{Taubes2000}). If $X$ is a symplectic manifold with $b_2^+(X)>1$ and $[A]$ is a second homology class of $X$ such that all of its representatives are connected and $g([A])\neq 1$, then the relation between the Gromov-Witten invariants and the Gromov invariants of Taubes is

\begin{center} $$Gr^X([A])=GW^X_{[A],k_{[A]}}(PD([point]),\cdots,PD([point]))$$

\vspace{1mm}

\end{center} where $k_{[A]}=[A]^2+1-g([A])$ and $PD([point])$ is repeated $k_{[A]}$ times. 

The number of points in $X$ for $[E]$ is determined by Taubes as $k_{[E]}=[E]^2+1-g([E])=0$. Keeping in mind that the representative for $[E]$ is connected, according to the formula which gives the relation between the Gromov invariants of Taubes and the Seiberg-Witten invariants (\cite{Fintushel2006a}), $GW_{[E],0}^X= Gr^X([E])=SW(2[E]+c_1(X))=\pm1$. The last equality is justified by the blowup formula for the Seiberg-Witten invariants ($c_1(X)=c_1(Y)-[E]$).}

The similarity of the calculations in Example~\ref{example4} and to the calculations in Lemma~\ref{theExample} brings to mind that there may be relation between the invariants of the underlying smooth structure of a $4$-manifold $X$ and the Gromov-Witten invariants of $X\times S^2$. Unfortunately, we see that no interesting relation may occur in the genus zero case.

\section{Exotic Symplectic Manifolds in Dimension Six} \label{exoticSix} 

Example~\ref{RuansExample} is brought first by Ruan in \cite{Ruan1994} in the context of exotic symplectic structures on smooth $6$-manifolds. See \cite{McDuff2004} (page 335) for another explanation of this example. The manifolds in these sources have $b_2^+=1$.

{\example \label{RuansExample} Let $X_1$ be $\CP\#_8\CPB$ and $X_2= B_8$ be the Barlow surface. $X_1$ and $X_2$ are homeomorphic (\cite{Freedman1982},\cite{Donaldson1983}) but they are not diffeomorphic (\cite{Kotschick1989}). The Barlow surface is minimal and $\CP\#_8\CPB$ is not minimal. So $X_1\times S^2$ and $X_2\times S^2$ are not symplectic deformation equivalent (\cite{Ruan1994}).} \\

The next lemma is based on this example.

{\lemma \label{theExample} Let $Y$ be a symplectic $4$-manifold and $X$ be its blowup, i.e $X=Y\#\CPB$. Let $[E]$ be the class of the exceptional sphere and let $\overline{[E]}$ denote the pushforward of the homology class $[E]$ under the inclusion map in $H_2(X\times S^2)$. Then $GW^{X\times S^2}_{\overline{[E]},1}(PD(\overline{[E]}))=-1$.}
\begin{proof} Let $J$ be a compatible almost complex structure on $X\times S^2$. The dimension of ${\M}^{X\times S^2}_{\overline{[E]},1}$ is $2c_1(X\times S^2)\overline{[E]}+2=2c_1(X)[E]+2$ which is equal to four, and the dimension condition is satisfied. By positivity of intersections of $J$-holomorphic curves in an almost complex manifold, there is only one curve which represents $[E]$ in $X$, which will be denoted by $E$. For each point of $S^2$ factor, we have the curve $\overline{E}$ in $X\times\cdot\,\,$, where $\overline{E}$ is the image of $E$ in $X\times S^2$. If we put the condition of passing through a marked point, this adds two real dimensions to the moduli space for the freedom of choosing a point on the sphere $E$. The moduli space is diffeomorphic to $\overline{E}\times S^2$ in $X\times S^2$ which is compact and smooth and the evaluation map is the diffeomorphism. In this case the intersection of the pseudocycles in the definition of Gromov-Witten invariants are in fact an intersection of cycles in $M$. Therefore $GW^{X\times S^2}_{\overline{[E]},1}(PD(\overline{[E]}))$ is equal to $[\overline{E}\times S^2]\cdot \overline{[E]}$ which is $-1$.
\end{proof}

When $k$ is zero, as discussed at the end of Section~\ref{GW}, the dimension condition does not hold since $c_1(X\times S^2)\overline{[E]}$ is nonzero. Thus $GW^{X\times S^2}_{\overline{[E]},0}$ is zero.

The following theorem extends Lemma~\ref{theExample} to the cases where $k$ is greater than one. See also \cite{Ruan2002}.

{\theorem \label{6pm1} Let $Y$ be a symplectic $4$-manifold, $X$ be $Y\#\CPB$ and $k\geq1$. Then $$GW^{X\times S^2}_{\overline{[E]},k}(PD(\overline{[E]}),PD([\overline{E}\times S^2]),\cdots,PD([\overline{E}\times S^2]))=(-1)^k$$}
\begin{proof} The dimension of ${\M}^{X\times S^2}_{\overline{[E]},k}$ is $2c_1(X\times S^2)\overline{[E]}+2k=2+2k$. If $k$ is one, then by Lemma~\ref{theExample} $GW^{X\times S^2}_{\overline{[E]},1}(PD(\overline{[E]}))$ is $-1$. $PD([\overline{E}\times S^2])$ is a degree two cohomology class of $X\times S^2$. So we can apply the divisor axiom (Lemma~\ref{divisorAxiom} ). Applying the divisor axiom inductively, we conclude that $GW^{X\times S^2}_{\overline{[E]},k}(PD(\overline{[E]}),PD(\overline{[E}\times S^2]),\cdots,PD(\overline{[E}\times S^2]))$ is $(-1)^k$.
\end{proof}

The second homology classes of $X\times S^2$ are pushforwards of the second homology classes of $X$ or $[\,\,\cdot\times S^2]$. The next theorem is on the former classes.

{\theorem \label{6A=E} Let $X$ be a symplectic $4$-manifold with $b_2^+(X)>1$ and $J$ be a generic almost complex structure on $X$ which is compatible with the symplectic structure. Let $[A]$ be a  nonzero second homology class of $X$ which can be represented by an embedded, connected $J$-holomorphic sphere. Let $\overline{[A]}$ be the pushforward of $[A]$ in $H_2(X\times S^2;\Z)$ and $\alpha_1,\cdots,\alpha_k$ be cohomology classes of $X\times S^2$. If $GW^{X\times S^2}_{\overline{[A]},k}(\alpha_1,\cdots,\alpha_k)$ is nonzero, then $[A]$ is the homology class of an exceptional sphere in $X$.}
\begin{proof} This result is a straightforward consequence of Theorem~\ref{4A=E}.
\end{proof}

The next theorem says that Ruan's example is the only meaningful application of genus zero Gromov-Witten invariants in the case of stabilized $4$-manifolds.

{\theorem \label{main} Let $X$ be a simply connected, symplectic $4$-manifold with $b_2^+(X)>1$ and $J$ be a generic almost complex structure on $X$ which is compatible with the symplectic structure. Let $[A]$ be a nonzero second homology class of $X$ which can be represented by an embedded, connected $J$-holomorphic sphere. Let $\overline{[A]}$ be the pushforward of $[A]$ in $H_2(X\times S^2;\Z)$ and $\alpha_1,\cdots,\alpha_k$ be cohomology classes of $X\times S^2$. If $GW^{X\times S^2}_{\overline{[A]},k}(\alpha_1,\cdots,\alpha_k)$ is nonzero, then the following conditions are satisfied.
\begin{enumerate} \item For an exceptional sphere $E$ in $X$, $[A]$ is the homology class $[E]$ of $E$, 
\item For a unique $j$, $\alpha_j$ is a fourth cohomology class of $X\times S^2$ which evaluates nonzero on $[\overline{E}\times S^2]\in X\times S^2$ where $\overline{E}$ is the image of $E$ in $X\times S^2$ under the inclusion map, 
\item For all $i$ which are not equal to $j$, $\alpha_i$ is a second cohomology class of $X\times S^2$ which evaluates nonzero on $\overline{[E]}\in X\times S^2$.
\end{enumerate}
}
\begin{proof} Assume that $GW^{X\times S^2}_{\overline{[A]},k}(\alpha_1,\cdots,\alpha_k)$ is not zero. Theorem~\ref{6A=E} imposes that $[A]$ should be the homology class of an exceptional sphere $E$ in $X$, i.e. $[A]$ is identical with $[E]$ in $H_2(X;\Z)$. By the dimension condition, the sum of degrees of $\alpha_i$ should be equal to the dimension of the moduli space ${\M}^{X\times S^2}_{\overline{[A]},k}$, which is $2n-6+2c_1(X\times S^2)\overline{[A]}+2k=2+2k$. $X\times S^2$ is simply connected, thus its odd cohomology groups are trivial, and all $\alpha_i$'s ($1\leq i\leq k$) have even degrees. By the fundamental class axiom (Lemma~\ref{fundamentalClassAxiom}), in order to get a nonzero invariant, there must be one class of degree four, and the remaining ones are of degree two. 

Without loss of generality, since there is no odd degree cohomology class, assume that $\alpha_1$ is the fourth degree class. 

The moduli space ${\M}^{X\times S^2}_{\overline{[E]},1}$ is diffeomorphic to $\overline{E}\times S^2$ in $X\times S^2$ which is compact and smooth and the evaluation map is the diffeomorphism. In this case the intersection of the pseudocycles in the definition of Gromov-Witten invariants are in fact an intersection of cycles in $M$. Therefore $GW^{X\times S^2}_{\overline{[E]},1}(\alpha_1)$ is equal to $\alpha_1\cdot[\overline{E}\times S^2]$, the evaluation of $\alpha_1$ on $[\overline{E}\times S^2]$, and is nonzero only if the latter is nonzero. 

Now let us turn to the cohomology classes $\alpha_i$, $2\leq i\leq k$. Each $\alpha_i$ is of degree two, so the divisor axiom is applicable and the result follows.

\end{proof}

{\corollary \label{GW6values} Let $k$ be a positive integer, $\alpha_1\in H^4(X;\Z)$ and $\alpha_i\in H^2(X;\Z)$ for all $i\in\Z$ such that $2\leq i\leq k$. If $X$ is a simply connected, symplectic $4$-manifold, then $$GW^{X\times S^2}_{\overline{[A]},k}(\alpha_1,\cdots,\alpha_k) = (\alpha_1\cdot[\overline{E}\times S^2])\cdot(\alpha_k\cdot \overline{[E]})\cdots(\alpha_k\cdot \overline{[E]})$$}

{\remark Let $X$ be a $4$-manifold as in Theorem~\ref{main}. If $GW^{X\times S^2}_{\overline{[A]},k}(\alpha_1,\cdots,\alpha_k)$ be nonzero for a nonzero second homology class $[A]$ of $X$ which can be represented by an embedded, connected $J$-holomorphic sphere and for some cohomology classes $\alpha_1,\cdots,\alpha_k$ of $X\times S^2$, then $X$ is a blowup of a $4$-manifold $Y$. $H_2(X;\Z)$ is isomorphic to the direct sum $H_2(Y;\Z)\oplus H_2(\CPB;\Z)$. By a slight abuse of notation, suppose that $H_2(\CPB;\Z)$ is generated by $[E]$. Then $H_2(X\times S^2;\Z)$ is isomorphic to the direct sum $H_2(Y;\Z)\oplus H_2(\CPB;\Z) \oplus H_2(S^2;\Z)$. A generator of $H_2(X\times S^2;\Z)$ is either the pushforward of a generator $[B]$ of $H_2(Y;\Z)$ under the inclusion map into $X\times S^2$, $\overline{[E]}$ or $[\,\,\cdot\times S^2]$. If the Poincar\'e dual of the degree four cohomology class $\alpha_1$ in the proof of the theorem is written as a linear combination of these generators, then the coefficient of $\overline{[E]}$ can not be zero. That is $\alpha_1$ has a $PD(\overline{[E]})$ term. A similar argument applies to $\alpha_i$ ($2\leq i\leq k$) and $[\overline{E}\times S^2]$.}\\

A corollary to this theorem is the existence of exotic symplectic deformation types on a fixed smooth $6$-manifold. 

{\corollary \label{deformationEq} Let $X_1$ and $X_2$ be homeomorphic symplectic $4$-manifolds with $b_2^+(X_1)>1$ and $b_2^+(X_2)>1$. If $X_1$ is not minimal and $X_2$ is minimal, then $X_1\times S^2$ and $X_2\times S^2$ are diffeomorphic but they are not symplectic deformation equivalent.}
\begin{proof} In Lemma~\ref{theExample}, one of the invariants of $X_1\times S^2$, $GW^{X\times S^2}_{\overline{[E]},1}(PD(\overline{[E]}))$, is found to be nonzero. 

$X_1\times S^2$ and $X_2\times S^2$ are diffeomorphic (\cite{Wall1966},\cite{Wall1967a},\cite{Jupp1973}). Let $h:X_1\times S^2\rightarrow X_2\times S^2$ be a diffeomorphism. The diffeomorphism $h$ induces an isomorphism oo the homology, the cohomology and the triple intersection forms. Under this isomorphism $c_1(X_1\times S^2)$ is taken to $c_1(X_2\times S^2)$, and by Theorem 9 of \cite{Wall1966} homotopy class of the compatible almost complex structures is preserved. Theorem\,\ref{main} implies that for a generic almost complex structure $X_2\times S^2$ has all its corresponding genus zero Gromov-Witten  invariants zero since $X_2$ is minimal. Therefore the symplectic structures on $X_1\times S^2$ and $X_2\times S^2$ are not symplectic deformation equivalent.
\end{proof}

As another outcome to Theorem~\ref{main} we see that for two minimal, simply connected, symplectic $4$-manifolds $X_1$ and $X_2$ such that $b_2^+(X_1)>1$ and $b_2^+(X_2)>1$, genus zero Gromov-Witten invariants for the pushforwards of second homology classes with minimal genera zero can not distinguish the symplectic structures on $X_1\times S^2$ and $X_2\times S^2$.

\bibliography{library}

\end{document}